\documentclass[12pt,reqno]{amsart} % amsart

\usepackage{mathrsfs,amsfonts,amsthm,amsmath,amssymb,thmtools} % thmtools
\usepackage[dvipsnames]{xcolor} % dvipsnames：支持68种颜色（CMYK）; svgnames：支持151种颜色（RGB）; xllnames：支持317种颜色（RGB）
\usepackage{enumitem}                % 精确控制列表环境 (如 itemize, enumerate, description) 样式和间距
\usepackage{lineno}                 % to add line numbers for correction

\usepackage[numbers,sort&compress]{natbib}  % 连续编号参考文献

\usepackage{graphics}
\usepackage{graphicx}   %从文件夹引用图片
\graphicspath{{images/}} %当前目录下图片所在文件夹
\usepackage{subfig}     %将图片插入到文档
\usepackage{rotating}

\usepackage[hyperindex, pdfencoding = auto, psdextra, bookmarksdepth = 4]{hyperref}  % 书签与链接的交叉引用
\hypersetup{
    bookmarksopen      = true,
    bookmarksopenlevel = 1,
    bookmarksnumbered  = true,
    pdftitle    = {},
    pdfauthor   = {},
    linktoc     = page,
    anchorcolor = green,
    colorlinks  = true,        % Colours links instead of ugly boxes,
    linkcolor   = black,
    citecolor   = blue,
    filecolor   = blue,
    urlcolor    = magenta,
}

\usepackage[capitalise, noabbrev, nameinlink]{cleveref}    % cleveref 宏包必须在 hyperref 宏包之后加载
\theoremstyle{plain}

\newtheorem{theorem}{Theorem}[section]
\newtheorem{lemma}[theorem]{Lemma}

\newtheorem{proposition}[theorem]{Proposition}
\newtheorem{corollary}[theorem]{Corollary}

\theoremstyle{definition}   % 来自amsthm宏包的定理样式设置命令 \theoremstyle{风格选项:  plain, definition, remark}
\makeatletter
    \def\subsection{\@startsection{subsection}{2}%
    \z@{.5\linespacing\@plus.7\linespacing}{.3\linespacing}%
    {\normalfont\bfseries}}
\makeatother % 修改subsection换行的命令

\makeatletter                                         % make之间的部分允许访问内部命令
    \newcommand{\LeftEqNo}{\let\veqno\@@leqno}        % 让公式编号出现在左侧
\makeatother

\numberwithin{equation}{section}                      % 让公式按 section 进行编号

\begin{document}

\

\vspace{-2cm}

\title[Trace-class perturbations of Voiculescu's free semicircular families]
{Trace-class perturbations of Voiculescu's free semicircular families}

\author{Minghui Ma}
\address{Minghui Ma, School of Mathematical Sciences, Dalian University of Technology, Dalian, 116024, China}
\email{minghuima@dlut.edu.cn}

\author{Rui Shi}
\address{Rui Shi, School of Mathematical Sciences, Dalian University of Technology, Dalian, 116024, China}
\email{ruishi@dlut.edu.cn, ruishi.math@gmail.com}

\author{Tianze Wang}
\address{Tianze Wang, School of Mathematical Sciences, Dalian University of Technology, Dalian, 116024, China}
\email{swan0108@mail.dlut.edu.cn}


\makeatletter
\@namedef{subjclassname@2020}{\textup{2020} Mathematics Subject Classification}
\makeatother

\begin{abstract}
Let $\{s_j\}_{j\in\Lambda}$ be the standard free semicircular family on the full Fock space $\mathcal H=\mathcal F(\ell^2(\Lambda))$, where $|\Lambda|\geqslant2$.
For any family $\{K_j\}_{j\in\Lambda}$ of self-adjoint trace-class operators with $\gamma=\sum_{j\in\Lambda}\|K_j\|_1<\infty$, we prove that the type $\mathrm{II}_1$ summand $\mathcal P$ of the von Neumann algebra $W^*(\{s_j+K_j\colon j\in\Lambda\})$ is isomorphic to the free group factor $\mathcal L(\mathbb F_\Lambda)$.
As an application, we show that the set of irreducible operators on a separable infinite-dimensional Hilbert space is not dense in the trace-class norm.
For each fixed $\gamma>0$, the possible values of the Murray--von Neumann dimension $\dim_{\mathcal P}(z_{\mathcal P}\mathcal H)$ form the interval $[1,1+\gamma)$, where $z_{\mathcal P}$ is the identity of $\mathcal P$.
\end{abstract}

\maketitle
%\end{frontmatter}

\section{Introduction}

The classical Weyl--von Neumann theorem \cite{Weyl-09,vN-35} states that every self-adjoint operator on a separable infinite-dimensional Hilbert space is diagonalizable up to an arbitrarily small Hilbert-Schmidt perturbation.
Kuroda \cite{Kuroda-58} extended the Weyl--von Neumann theorem to every unitarily invariant norm not equivalent to the trace-class norm, including Schatten $p$-norms with $1<p<\infty$.
At $p=1$, the Kato--Rosenblum theorem \cite{Kato-57,Rosenblum-57} states that the absolutely continuous parts of two self-adjoint operators with trace-class difference are unitarily equivalent.
In particular, a nonzero absolutely continuous part obstructs diagonalization by a trace-class perturbation.
For bounded self-adjoint $n$-tuples with commuting entries, $n\geqslant 2$, Voiculescu \cite{Voiculescu-79,Voiculescu-81} proved that coordinatewise perturbations in the Lorentz ideal $\mathcal C_n^-$ preserve the joint Lebesgue absolutely continuous part up to unitary equivalence, provided the perturbed tuple also has commuting entries. He \cite[Theorem 3.1]{Voiculescu-20} later proved the same conclusion for coordinatewise perturbations in $\mathcal C_{p_j}^-$, where $p_j>1$ and $\sum_{j=1}^n p_j^{-1}=1$. His recent survey \cite{Voiculescu-25}  discusses wave operators for commuting tuples and formulates two questions for pairs. For $n\geqslant 3$, Bhutani and Voiculescu \cite{BV-26} proved the existence of a strong wave operator on the joint absolutely continuous subspace for possibly unbounded, strongly commuting tuples under coordinatewise $\mathcal C_n^-$ perturbations.

% For $n\geqslant2$, Voiculescu \cite{Voiculescu-79,Voiculescu-81} proved the corresponding result for commuting $n$-tuples of self-adjoint operators when the perturbation in each coordinate belongs to the Lorentz ideal $\mathcal C_n^-$.

Free probability theory was introduced by Voiculescu \cite{VDN-92} in his study of free group factors.
In the free central limit theorem \cite[Theorem~3.5.1, p.~29]{VDN-92}, the semicircle law plays the same role as the Gaussian law in the classical central limit theorem.
Voiculescu's construction on the full Fock space gives a free semicircular family whose generated von Neumann algebra is isomorphic to the free group factor on the same index set \cite[Theorem~2.6.2, p.~17]{VDN-92}.
The Kato--Rosenblum theorem applies to each semicircular generator separately.
Inspired by Voiculescu's work \cite{Voiculescu-79,Voiculescu-81}, we study simultaneous trace-class perturbations of a free semicircular family.

Let $\Lambda$ be a set with $|\Lambda|\geqslant2$, and consider the full Fock space
\begin{equation*}
    \mathcal H
    =\mathcal F(\ell^2(\Lambda))
    =\mathbb C\Omega\oplus\bigoplus_{n=1}^{\infty}\ell^2(\Lambda)^{\otimes n}.
\end{equation*}
For each $j\in\Lambda$, let $l_j$ and $r_j$ be the \emph{left and right creation operators} associated with the $j$-th standard unit vector, respectively.
Their adjoints $l_j^*$ and $r_j^*$ are the corresponding \emph{annihilation operators}.
On a Fock word, $l_j^*$ removes the first letter if it is $j$, and gives zero otherwise; $r_j^*$ acts similarly on the last letter.
Both annihilate the vacuum vector, that is, $l_j^*\Omega=r_j^*\Omega=0$.

Let $s_j=l_j+l_j^*$.
The operators $\{s_j\}_{j\in\Lambda}$ form a standard free semicircular family.
By \cite[Example~3.4.2, p.~26]{VDN-92}, the distribution of each $s_j$ with respect to the vacuum state $T\mapsto\langle T\Omega,\Omega\rangle$ is the standard semicircle law $\mu$, given by
\begin{equation*}
    d\mu(t)=\frac{1}{2\pi}\sqrt{4-t^2}\,dt,
    \qquad -2\leqslant t\leqslant2.
\end{equation*}
In particular, $\tau(s_j^2)=1$.
The vacuum state restricts to a faithful normal tracial state $\tau$ on $W^*(\{s_j\colon j\in\Lambda\})$.
The standard identifications
\begin{equation}\label{equ identification}
    W^*(\{s_j\colon j\in\Lambda\})
    \cong\mathcal L(\mathbb F_\Lambda),
    \qquad
    \mathcal H\cong L^2(\mathcal L(\mathbb F_\Lambda),\tau)
\end{equation}
will be used throughout, where $\mathbb F_\Lambda$ is the free group on $\Lambda$ and $\tau$ corresponds to its canonical group trace; see \cite[Theorem~2.6.2, p.~17]{VDN-92}.
Moreover, the inner product is given by $\langle x,y\rangle=\tau(y^*x)$ for $x,y\in\mathcal L(\mathbb F_\Lambda)$.

For a trace-class operator $K$ on $\mathcal H$, we write $\|K\|_1=\operatorname{Tr}(|K|)$.
Our first theorem shows that when the semicircular generators are perturbed by self-adjoint operators with summable trace norms, the resulting von Neumann algebra has a direct summand isomorphic to the original free group factor.

\begin{theorem}\label{thm main}
    Let $\{s_j\}_{j\in\Lambda}$ be the standard free semicircular family defined above.
    Suppose that $\{K_j\}_{j\in\Lambda}$ is a family of self-adjoint trace-class operators on $\mathcal H$ satisfying
    \begin{equation}\label{equ sum-K}
        \gamma=\sum_{j\in\Lambda}\|K_j\|_1<\infty.
    \end{equation}
    Let $S_j=s_j+K_j$ for $j\in\Lambda$ and $\mathcal N=W^*(\{S_j\colon j\in\Lambda\})$.
    Then $\mathcal N$ has a central decomposition
    \begin{equation*}
        \mathcal N=\mathcal P\oplus\mathcal Q\oplus\mathcal D,
    \end{equation*}
    where $\mathcal P=z_{\mathcal P}\mathcal N$ is isomorphic to $\mathcal L(\mathbb F_\Lambda)$ for a nonzero central projection $z_{\mathcal P}$ in $\mathcal N$, $\mathcal Q$ is either zero or properly infinite, and $\mathcal D$ acts on a finite-dimensional reducing subspace of $\mathcal H$.
    Moreover, there exists an isomorphism $\theta\colon\mathcal P\to\mathcal L(\mathbb F_\Lambda)$ such that
    \begin{equation*}
        \theta(z_{\mathcal P}S_j)=s_j\qquad(j\in\Lambda).
    \end{equation*}
    If $\Lambda$ is infinite, then $\mathcal D=0$.
\end{theorem}

When $\Lambda$ is finite, the summability condition \eqref{equ sum-K} is automatic.
When $\Lambda$ is uncountable, it implies that $K_j=0$ except for countably many $j\in\Lambda$.

A long-standing problem in operator theory concerns the density of irreducible operators on separable infinite-dimensional Hilbert spaces in the trace-class norm.
Halmos \cite{Halmos-68} proved that irreducible operators form a dense $G_\delta$ subset in the operator norm topology.
In the last paragraph of \cite[Section 1]{Halmos-68}, Douglas observed the density of irreducible operators in the Hilbert-Schmidt norm by virtue of the Weyl--von Neumann theorem.
Herrero proved in \cite[Lemma 4.33]{Herrero-82-book} that the set of irreducible operators is dense in the Schatten $p$-norm for every $p>1$.
The authors in \cite[Theorem~1.2]{FJMSSW-25} reduced the long-standing density problem to the case of generators of cyclic type $\mathrm{II}_1$ factors.
Applying \Cref{thm main}, we answer this problem negatively in \Cref{cor irreducible}, where we establish a connection between operator theory and von Neumann algebras.

The unitary equivalence in the Kato--Rosenblum theorem implies that the absolutely continuous spectral multiplicity functions of self-adjoint operators are invariant under trace-class perturbation.
For commuting tuples, Voiculescu \cite{Voiculescu-18} also related the integral of the joint absolutely continuous spectral multiplicity function to the quasicentral modulus for a family of Lorentz ideals.
In the setting of \Cref{thm main}, we may measure the multiplicity of the representation of $\mathcal P$ on $z_{\mathcal P}\mathcal H$ by the Murray--von Neumann dimension $\dim_{\mathcal P}(z_{\mathcal P}\mathcal H)$, as in \cite[Chapter~V, Section~10$\alpha$]{Connes-94-book}.
For the unperturbed family, this dimension is $1$, since $\mathcal H\cong L^2(\mathcal L(\mathbb F_\Lambda),\tau)$ as a left $\mathcal L(\mathbb F_\Lambda)$-module.
The next theorem shows that this dimension can increase under trace-class perturbations and determines possible values of this dimension for each fixed total trace norm $\gamma=\sum_{j\in\Lambda}\|K_j\|_1$.

\begin{theorem}\label{thm multiplicity-range}
    With the notation of \Cref{thm main}, let $d=\dim_{\mathcal P}(z_{\mathcal P}\mathcal H)$.
    If $\gamma=0$, then $d=1$.
    For each $\gamma>0$, as $\{K_j\}_{j\in\Lambda}$ varies over all self-adjoint trace-class families satisfying $\sum_{j\in\Lambda}\|K_j\|_1=\gamma$, the dimension $d$ takes exactly the values in $[1,1+\gamma)$.
\end{theorem}

In \Cref{sec pre}, we record the commutation relations on the full Fock space and construct a tracial functional used to bound the Murray--von Neumann dimension.
In \Cref{sec compact}, we show that, for compact perturbations, a nonzero intertwiner between the original and perturbed families yields a central summand isomorphic to $\mathcal L(\mathbb F_\Lambda)$.
In \Cref{sec trace-class}, we construct such an intertwiner by finite-rank approximation under the summability assumption \eqref{equ sum-K}. We then prove \Cref{thm main} and apply it to the density problem for irreducible operators.
Finally, in \Cref{sec multiplicity}, we prove \Cref{thm multiplicity-range} by establishing the dimension bounds and, for each $\gamma>0$, realizing every value in $[1,1+\gamma)$ by a positive finite-rank perturbation of a single generator.

\section{Preliminaries}
\label{sec pre}

\subsection{Commutant and commutation identities}

We record the commutant and commutation identities used below.
Under the preceding identifications \eqref{equ identification}, \cite[Remark~2.6.7, p.~20]{VDN-92} gives
\begin{equation*}
    \mathcal L(\mathbb F_\Lambda)'
    =W^*(\{r_j+r_j^*\colon j\in\Lambda\})
    =R(\mathcal L(\mathbb F_\Lambda)).
\end{equation*}
Here $R_x\xi=\xi x$, and the canonical trace $\tau'$ on the commutant satisfies $\tau'(R_x)=\tau(x)$.
Let $P_\Omega$ be the rank-one projection onto $\mathbb C\Omega$, and let $e_j=l_j\Omega=r_j\Omega$ for $j\in\Lambda$.
Using the definitions of the left and right creation operators in \cite[Example~1.5.8, p.~6, and Remark~2.6.7, p.~20]{VDN-92}, a direct computation on Fock words gives
\begin{equation*}
    [l_j,r_i]=0,
    \qquad
    [l_j^*,r_i]=\delta_{ij}P_\Omega,
    \qquad
    [l_j,r_i^*]=-\delta_{ij}P_\Omega\quad(i,j\in\Lambda).
\end{equation*}

\subsection{Murray-von Neumann dimension}

We record a well-known result on the Murray-von Neumann dimension.

\begin{lemma}\label{lem cyclic-dimension}
Let $(\mathcal M,\tau)$ be a type $\mathrm{II}_1$ factor and let $\mathcal H$ be a normal left $\mathcal M$-module.
Suppose that $\xi\in\mathcal H$ and $p\in\mathcal M$ is a projection
such that $p\xi=\xi$. Then
\begin{equation*}
    \dim_{\mathcal M}\overline{\mathcal M\xi}\leqslant\tau(p).
\end{equation*}
In particular, a submodule generated by $n$ vectors has $\mathcal M$-dimension at most $n$.
\end{lemma}

\begin{proof}
Let $g\in L^1(\mathcal M,\tau)_+$ be defined by
\begin{equation*}
    \tau(xg)=\langle x\xi,\xi\rangle
    \qquad(x\in\mathcal M).
\end{equation*}
Then $\tau((1-p)g)=\|(1-p)\xi\|^2=0$, and hence $g=pgp$.
Since
\begin{equation*}
	\|x\xi\|^2=\tau(x^*xg)=\|xg^{1/2}\|_2^2,
\end{equation*}
the map $x\xi\mapsto xg^{1/2}$ extends to an isometric $\mathcal M$-module embedding of $\overline{\mathcal M\xi}$ into $L^2(\mathcal M,\tau)p$.
The latter module has dimension $\tau(p)$.
Taking $p=1$, every cyclic submodule has dimension at most $1$.
The last assertion follows from subadditivity of dimension for sums of submodules.
\end{proof}

\subsection{A tracial functional}
\label{subsec tracial}

In this subsection, $\mathcal H$ denotes the full Fock space $\mathcal F(\ell^2(\Lambda))$.
Let $\mathcal E$ denote the set of all operators $X\in \mathcal B(\mathcal H)$ such that
\begin{equation*}
	[s_j,X]\in\mathcal S_1(\mathcal H)\quad(j\in\Lambda),
	\qquad
	\sum_{j\in\Lambda}\|[s_j,X]\|_1<\infty.
\end{equation*}
It is clear that $\mathcal E$ is a unital $*$-subalgebra of $\mathcal B(\mathcal H)$.
For a finite tuple, this is the commutant modulo the trace class considered in \cite[Section~2.1]{Voiculescu-19}.
Motivated by Voiculescu's trace pairing \cite[Proposition~3.1]{Voiculescu-17}, we define a linear functional by
\begin{equation*}
	\varphi(X)=\langle X\Omega,\Omega\rangle+\sum_{j\in\Lambda}\operatorname{Tr}([s_j,X]\,l_j)
	\qquad(X\in\mathcal E).
\end{equation*}
For each $X\in\mathcal E$, at most countably many commutators $[s_j,X]$ are nonzero, and the sum converges absolutely since $\|l_j\|=1$.
It is clear that $\varphi$ agrees with the vacuum trace $\tau'$ on $\mathcal L(\mathbb F_\Lambda)'$.
Moreover, if $\Lambda$ is finite, then $P_\Omega\in\mathcal E$ and $[s_j,P_\Omega]l_j=-P_\Omega$ for every $j\in\Lambda$.
Thus, $\varphi(P_\Omega)=1-|\Lambda|<0$.

We first verify the identity $\varphi(AX)=\varphi(XA)$ for $X\in\mathcal E$ when $A$ is a right creation operator or its adjoint.

\begin{lemma}\label{lem r-in-F}
Let $\mathcal F=\{A\in\mathcal E\colon\varphi([A,X])=0~\text{for every}~X\in\mathcal E\}$.
Then $\mathcal F$ is a unital subalgebra of $\mathcal E$, and $r_j,r_j^*\in\mathcal F$ for every $j\in\Lambda$.
\end{lemma}

\begin{proof}
If $A,B\in\mathcal F$ and $X\in\mathcal E$, then $BX,XA\in\mathcal E$, and
\begin{equation*}
	\varphi([AB,X])
	=\varphi([A,BX])+\varphi([B,XA])=0.
\end{equation*}
Thus, $\mathcal F$ is a unital subalgebra of $\mathcal E$.

Let $\omega_\Omega(T)=\langle T\Omega,\Omega\rangle$ for $T\in\mathcal B(\mathcal H)$.
For $j\in\Lambda$ and $A,X\in\mathcal E$, we have
\begin{equation*}
	[s_j,[A,X]]=[s_j,A]X-X[s_j,A]+A[s_j,X]-[s_j,X]A,
\end{equation*}
and hence
\begin{equation*}
	\varphi([A,X])=\omega_\Omega([A,X])
	+\sum_{j\in\Lambda}\operatorname{Tr}\bigl([s_j,A][X,l_j]+[s_j,X][l_j,A]\bigr).
\end{equation*}
Recall that $[l_j,r_i]=0$, $[l_j^*,r_i]=\delta_{ij}P_\Omega$, and $[l_j,r_i^*]=-\delta_{ij}P_\Omega$.
Since $l_j\Omega=r_j\Omega=e_j$ and $l_j^*\Omega=r_j^*\Omega=0$, we obtain
\begin{align*}
	\varphi([r_i,X])
	&=\omega_\Omega([r_i,X])+\operatorname{Tr}(P_\Omega[X,l_i])=0,\\
	\varphi([r_i^*,X])
	&=\omega_\Omega([r_i^*,X])
      -\operatorname{Tr}(P_\Omega[X,l_i])
      -\operatorname{Tr}([s_i,X]P_\Omega)=0.
\end{align*}
Thus, $r_i,r_i^*\in\mathcal F$.
\end{proof}

Next we show that $\varphi$ is tracial on $\mathcal E$.
For $\eta_1,\eta_2\in\mathcal H$, we write $\eta_1\,\hat{\otimes}\,\eta_2$ for the operator $\xi\mapsto\langle\xi,\eta_2\rangle\eta_1$.

\begin{proposition}\label{prop full-trace}
The linear functional $\varphi$ is tracial on $\mathcal E$, that is,
\begin{equation*}
	\varphi(AX)=\varphi(XA)\qquad(A,X\in\mathcal E).
\end{equation*}
\end{proposition}

\begin{proof}
We need to prove $\mathcal E=\mathcal F$, where $\mathcal F$ is given by \Cref{lem r-in-F}.
For each $X\in\mathcal E$, we have
\begin{equation*}
	\varphi([A,X])=\omega_\Omega([A,X])+\operatorname{Tr}(D_XA)
	+\sum_{j\in\Lambda}\operatorname{Tr}([s_j,A][X,l_j])
	\qquad(A\in\mathcal E),
\end{equation*}
where $D_X=\sum_{j\in\Lambda}[[s_j,X],l_j]\in\mathcal S_1(\mathcal H)$ and $\|D_X\|_1\leqslant 2\sum_{j\in\Lambda}\|[s_j,X]\|_1$.
The first two terms are ultraweakly continuous in $A$ and the last term can be dominated by $2\|X\|\sum_j\|[s_j,A]\|_1$.
Therefore, $\mathcal F$ has the \emph{closure property}:
if a net $A_\lambda\in\mathcal F$ converges ultraweakly to $A\in\mathcal E$ and satisfies $\sum_{j\in\Lambda}\|[s_j,A_\lambda-A]\|_1\to0$, then $A\in\mathcal F$.

By \Cref{lem r-in-F}, the unital $*$-algebra generated by $\{r_j+r_j^*\colon j\in\Lambda\}$ lies in $\mathcal F$.
This algebra is ultraweakly dense in $\mathcal L(\mathbb F_\Lambda)'$.
Since every operator in $\mathcal L(\mathbb F_\Lambda)'$ commutes with each $s_j$, the closure property gives $\mathcal L(\mathbb F_\Lambda)'\subseteq\mathcal F$.

Fix $A\in\mathcal E$.
Since $\sum_{j\in\Lambda}\|[s_j,A]\|_1<\infty$, there are sequences $\{\xi_m\}_{m=1}^\infty$ and $\{\eta_m\}$ in $\mathcal H=L^2(\mathcal L(\mathbb F_\Lambda),\tau)$ and labels $i(m)\in\Lambda$ such that
\begin{equation*}
	[s_j,A]=\sum_{i(m)=j}\xi_m\,\hat{\otimes}\,\eta_m\quad(j\in\Lambda),
	\qquad
	\sum_{m=1}^\infty(\|\xi_m\|^2+\|\eta_m\|^2)<\infty.
\end{equation*}
Let $d=\sum_{m=1}^\infty(\xi_m^*\xi_m+\eta_m^*\eta_m)\in L^1(\mathcal L(\mathbb F_\Lambda),\tau)_+$.
For $t>0$, let
\begin{equation*}
	p_t=\chi_{[0,t]}(d)\in\mathcal L(\mathbb F_\Lambda).
\end{equation*}
Then $\xi_mp_t,\eta_mp_t\in\mathcal L(\mathbb F_\Lambda)$.
It follows that $R_{\xi_mp_t},R_{\eta_mp_t}^*\in\mathcal L(\mathbb F_\Lambda)'\subseteq\mathcal F$.
For $n\geqslant 1$, let
\begin{equation*}
	B_{t,n}=\sum_{m=1}^n R_{\xi_mp_t}\,r_{i(m)}\,R_{\eta_mp_t}^*\in\mathcal F.
\end{equation*}
The Cauchy--Schwarz inequality gives
\begin{equation*}
	\|B_{t,n}\|\leqslant
	\left\|p_t\left(\sum_{m=1}^n\xi_m^*\xi_m\right)p_t\right\|^{1/2}
	\left\|p_t\left(\sum_{m=1}^n\eta_m^*\eta_m\right)p_t\right\|^{1/2}
	\leqslant t.
\end{equation*}
The identity $[s_j,r_i]=\delta_{ij}P_\Omega$ gives
\begin{equation*}
	\sum_{j\in\Lambda}\|[s_j,B_{t,n}]-R_{p_t}[s_j,A]R_{p_t}\|_1
	\leqslant\sum_{m>n}\|\xi_mp_t\|\|\eta_mp_t\|\to0.
\end{equation*}
Let $B_t$ be a weak-operator cluster point of the bounded sequence $\{B_{t,n}\}_{n=1}^\infty$.
The preceding estimate gives $[s_j,B_t]=R_{p_t}[s_j,A]R_{p_t}$.
Since $B_{t,n}\in\mathcal F$, the closure property gives $B_t\in\mathcal F$.
Moreover, we have
\begin{equation*}
    [s_j,R_{p_t}AR_{p_t}-B_t]
    =R_{p_t}[s_j,A]R_{p_t}-[s_j,B_t]=0
    \qquad(j\in\Lambda).
\end{equation*}
Thus, $R_{p_t}AR_{p_t}-B_t\in\mathcal L(\mathbb F_\Lambda)'\subseteq\mathcal F$, and hence $R_{p_t}AR_{p_t}\in\mathcal F$.
Note that $R_{p_t}\to I$ in the strong-operator topology and $R_{p_t}AR_{p_t}\to A$ in the weak-operator topology.
Moreover, we have
\begin{equation*}
	\sum_{j\in\Lambda}\|[s_j,R_{p_t}AR_{p_t}-A]\|_1
	=\sum_{j\in\Lambda}\|R_{p_t}[s_j,A]R_{p_t}-[s_j,A]\|_1\to0\quad(t\to\infty).
\end{equation*}
Applying the closure property once more yields $A\in\mathcal F$.
\end{proof}

\section{Compact perturbations}
\label{sec compact}

We study compact perturbations in the standard representation of a general type $\mathrm{II}_1$ factor.
For a unital $C^*$-algebra $\mathcal A\subseteq\mathcal B(\mathcal H)$, a state $\rho$ on $\mathcal B(\mathcal H)$ is called \emph{$\mathcal A$-central} if $\rho(aT)=\rho(Ta)$ for every $a\in\mathcal A$ and $T\in\mathcal B(\mathcal H)$.

\begin{lemma}\label{lem compact-perturbation}
Let $\mathcal M$ be a type $\mathrm{II}_1$ factor acting on $\mathcal H=L^2(\mathcal M,\tau)$, let $\{x_j\}_{j\in\Lambda}$ be a self-adjoint generating family of $\mathcal M$, and let $\mathcal A$ be the unital $C^*$-algebra generated by this family.
Suppose that $\mathcal A$ has the unique tracial state $\tau|_{\mathcal A}$ and $\mathcal B(\mathcal H)$ admits no $\mathcal A$-central state.
For self-adjoint compact operators $\{K_j\}_{j\in\Lambda}$ on $\mathcal H$, let
\begin{equation*}
    X_j=x_j+K_j\quad(j\in\Lambda),
    \qquad
    \mathcal N=W^*(\{X_j\colon j\in\Lambda\}).
\end{equation*}
Then $\mathcal N$ has a central decomposition
\begin{equation*}
    \mathcal N=\mathcal P\oplus\mathcal Q\oplus\mathcal D,
\end{equation*}
where $\mathcal P$ is either zero or isomorphic to $\mathcal M$, $\mathcal Q$ is either zero or properly infinite, and $\mathcal D$ acts on a finite-dimensional reducing subspace of $\mathcal H$.
If $\mathcal P=z_{\mathcal P}\mathcal N\neq0$, then there is an isomorphism $\theta\colon\mathcal P\to\mathcal M$ satisfying $\theta(z_{\mathcal P}X_j)=x_j$ for every $j\in\Lambda$.
\end{lemma}

\begin{proof}
Let $\mathcal B$ be the unital $C^*$-algebra generated by $\{X_j\}_{j\in\Lambda}$, let $\mathcal I=\mathcal B\cap\mathcal K(\mathcal H)$, and let $\pi$ be the Calkin quotient.
Since $\pi(X_j)=\pi(x_j)$ for every $j\in\Lambda$, it follows that $\pi(\mathcal B)=\pi(\mathcal A)$.
The restriction $\pi|_{\mathcal A}$ is injective because $\mathcal A\subseteq\mathcal M$ and $\mathcal M\cap\mathcal K(\mathcal H)=\{0\}$.
Thus, $\eta=(\pi|_{\mathcal A})^{-1}\circ(\pi|_{\mathcal B})$ gives the short exact sequence
\begin{equation*}
	0\longrightarrow\mathcal I\longrightarrow\mathcal B
	\mathrel{\mathop{\longrightarrow}\limits^{\eta}}\mathcal A\longrightarrow 0,
	\qquad
	\eta(X_j)=x_j\quad(j\in\Lambda).
\end{equation*}
Uniqueness of the tracial state on $\mathcal A$ implies that $\tau_{\mathcal B}=\tau\circ\eta$ is the unique tracial state on $\mathcal B$ annihilating $\mathcal I$.

Since $\mathcal I$ is a two-sided ideal of $\mathcal B$, the weak-operator closure $\overline{\mathcal I}^{\mathrm{WOT}}$ is a two-sided ideal of $\mathcal N$.
Hence $\overline{\mathcal I}^{\mathrm{WOT}}=z_{\mathcal I}\mathcal N$ for a central projection $z_{\mathcal I}\in\mathcal N$.
Let $z_f$ be the largest central projection of $\mathcal N$ for which $z_f\mathcal N$ is finite, and define
\begin{equation*}
    z_{\mathcal P}=z_f(1-z_{\mathcal I}),
    \qquad
    z_{\mathcal Q}=1-z_f,
    \qquad
    z_{\mathcal D}=z_fz_{\mathcal I}.
\end{equation*}
Let $\mathcal P=z_{\mathcal P}\mathcal N$, $\mathcal Q=z_{\mathcal Q}\mathcal N$, and $\mathcal D=z_{\mathcal D}\mathcal N$.
These projections are mutually orthogonal and sum to $1$.
By the choice of $z_f$, $\mathcal Q$ is either zero or properly infinite.

Suppose that $z_{\mathcal P}\neq0$.
For every normal tracial state $\varphi$ on $\mathcal P$, the tracial state $b\mapsto\varphi(z_{\mathcal P}b)$ on $\mathcal B$ annihilates $\mathcal I$ and hence equals $\tau_{\mathcal B}$.
The ultraweak density of $z_{\mathcal P}\mathcal B$ in $\mathcal P=z_{\mathcal P}\mathcal N$ shows that $\mathcal P$ has at most one normal tracial state.
It follows that $\mathcal P$ is a finite factor and admits a unique normal faithful tracial state $\tau_{\mathcal P}$.
For $b\in\mathcal B$, we have
\begin{equation*}
    \tau_{\mathcal P}\bigl((z_{\mathcal P}b)^*(z_{\mathcal P}b)\bigr)
    =\tau_{\mathcal B}(b^*b)
    =\tau(\eta(b)^*\eta(b)).
\end{equation*}
By faithfulness of $\tau_{\mathcal P}$ and $\tau$, $z_{\mathcal P}b=0$ if and only if $\eta(b)=0$.
Hence $z_{\mathcal P}b\mapsto\eta(b)$ defines a trace-preserving $*$-isomorphism from $z_{\mathcal P}\mathcal B$ onto $\mathcal A$.
The induced unitary operator between the $L^2$-spaces identifies the left multiplication representations and extends this map to an isomorphism $\theta\colon\mathcal P\to\mathcal M$ with $\theta(z_{\mathcal P}X_j)=x_j$ for every $j\in\Lambda$.

Let $w$ be a nonzero central projection in $z_{\mathcal D}\mathcal N$.
Since $w\overline{\mathcal I}^{\mathrm{WOT}}=w\mathcal N$, we have $w\mathcal I\neq\{0\}$.
The absolute value of a nonzero operator in $w\mathcal I$ has a nonzero finite-rank spectral projection $q\in w\mathcal N$.
Choose a minimal projection $p$ in the finite-dimensional algebra $q\mathcal Nq$.
Then the central support $z$ of $p$ in $z_{\mathcal D}\mathcal N$ satisfies $z\leqslant w$, and $z\mathcal N$ is a finite type $\mathrm{I}$ factor.
The projection $z$ is a finite sum of mutually orthogonal projections equivalent to $p$, and hence has finite rank.

We claim that $z_{\mathcal D}\mathcal H$ is finite dimensional.
Otherwise, applying the preceding argument successively to the nonzero central projections $z_{\mathcal D}-\sum_{k<n}z_k$, we obtain an orthogonal sequence $\{z_n\}_{n=1}^\infty$ of nonzero finite-rank central projections in $z_{\mathcal D}\mathcal N$.
Define $\mathcal N$-central states $\{\rho_n\}_{n=1}^\infty$ by
\begin{equation*}
    \rho_n(T)=
    \frac{\operatorname{Tr}(z_nTz_n)}{\operatorname{Tr}(z_n)}
    \qquad(T\in\mathcal B(\mathcal H)).
\end{equation*}
Since $z_n\to 0$ in the strong-operator topology, we have $|\rho_n(K)|\leqslant\|z_nKz_n\|\to 0$ for every compact operator $K$.
Let $\rho_\infty$ be a weak-* cluster point of $\{\rho_n\}_{n=1}^\infty$.
Then $\rho_\infty$ is an $\mathcal N$-central state annihilating compact operators.
For every $j\in\Lambda$ and $T\in\mathcal B(\mathcal H)$, we have
\begin{equation*}
	\rho_\infty([x_j,T])=\rho_\infty([X_j,T])-\rho_\infty([K_j,T])=0.
\end{equation*}
This makes $\rho_\infty$ an $\mathcal A$-central state, a contradiction.
Hence $\dim(z_{\mathcal D}\mathcal H)<\infty$.
\end{proof}

The next lemma gives a criterion for the existence of a central summand isomorphic to the original factor.

\begin{lemma}\label{lem nonzero-intertwiner}
Let $\mathcal M$ and $\mathcal N$ be von Neumann algebras on $\mathcal H$ generated by self-adjoint families $\{x_j\}_{j\in\Lambda}$ and $\{X_j\}_{j\in\Lambda}$, respectively.
Suppose that $\mathcal M$ is a type $\mathrm{II}_1$ factor and there exists a nonzero $Y\in\mathcal B(\mathcal H)$ satisfying $X_jY=Yx_j$ for every $j\in\Lambda$.
Then there are a nonzero central projection $z\in\mathcal N$ and an isomorphism $\theta\colon z\mathcal N\to\mathcal M$ such that $\theta(zX_j)=x_j$ for every $j\in\Lambda$.
\end{lemma}

\begin{proof}
Let $Y=V|Y|$ be the polar decomposition, and let $P=V^*V$ and $Q=VV^*$.
The intertwining relations and their adjoints give
$Y^*Y\in\mathcal M'$ and $YY^*\in\mathcal N'$.
It follows that $P\in\mathcal M'$ and $Q\in\mathcal N'$.
Since $|Y|\in\mathcal M'$, we have $(X_jV-Vx_j)|Y|=0$, and hence $(X_jV-Vx_j)P=0$.
Thus, we have
\begin{equation*}
	(X_jQ)V=V(x_jP),\qquad j\in\Lambda.
\end{equation*}
Since $\mathcal M$ is a factor and $P\neq0$, its normal compression to $P\mathcal H$ is faithful.
The normal compression of $\mathcal N$ to $Q\mathcal H$ has kernel $(1-z)\mathcal N$ for a nonzero central projection $z$.
Conjugation by the unitary operator $V\colon P\mathcal H\to Q\mathcal H$ identifies the compressed algebras $\mathcal MP$ and $\mathcal NQ$ and gives the required isomorphism.
\end{proof}

We now apply the preceding lemmas to the standard free semicircular family.

\begin{proposition}\label{prop compact-perturbation}
For self-adjoint compact operators $\{K_j\}_{j\in\Lambda}$ on $\mathcal H=L^2(\mathcal L(\mathbb F_\Lambda),\tau)$, let
\begin{equation*}
    S_j=s_j+K_j\quad(j\in\Lambda),
    \qquad
    \mathcal N=W^*(\{S_j\colon j\in\Lambda\}).
\end{equation*}
If there exists a nonzero $Y\in\mathcal B(\mathcal H)$ such that $S_jY=Ys_j$ for every $j\in\Lambda$, then $\mathcal N$ has a central decomposition
\begin{equation*}
    \mathcal N=\mathcal P\oplus\mathcal Q\oplus\mathcal D,
\end{equation*}
where $\mathcal P=z_{\mathcal P}\mathcal N$ is isomorphic to $\mathcal L(\mathbb F_\Lambda)$, $\mathcal Q$ is either zero or properly infinite, and $\mathcal D$ acts on a finite-dimensional Hilbert space.
Moreover, the isomorphism $\theta\colon\mathcal P\to\mathcal L(\mathbb F_\Lambda)$ can be chosen with $\theta(z_{\mathcal P}S_j)=s_j$ for all $j\in\Lambda$.
\end{proposition}

\begin{proof}
Let $\mathcal A$ be the unital $C^*$-algebra generated by $\{s_j\}_{j\in\Lambda}$.
Recall that the standard semicircle law $\mu$ is given by
\begin{equation*}
    d\mu(t)=\frac{1}{2\pi}\sqrt{4-t^2}\,dt,
    \qquad -2\leqslant t\leqslant2.
\end{equation*}
Define $F(t)=\mu([-2,t])$ for $t\in[-2,2]$, and let
\begin{equation*}
	u_j=e^{2\pi\mathbf{i}F(s_j)}\in C^*(1,s_j)
	\qquad(j\in\Lambda).
\end{equation*}
Then $\{u_j\}_{j\in\Lambda}$ is a family of free Haar unitary operators and $W^*(u_j)=W^*(s_j)$ for every $j\in\Lambda$; see \cite[Lemma~2.6.5, p.~18]{VDN-92}.
Freeness identifies $\mathcal A$ with the reduced free product of $C^*(1,s_j)$ and the unital $C^*$-algebra generated by the remaining semicircular operators, with respect to their restricted traces.
By \cite[Proposition~3.1]{Avitzour-82}, $\tau|_{\mathcal A}$ is the unique tracial state on $\mathcal A$.
Moreover, an $\mathcal A$-central state on $\mathcal B(\mathcal H)$ would make the canonical trace on the reduced free group $C^*$-algebra $C_r^*(\mathbb F_\Lambda)$ amenable, where $C_r^*(\mathbb F_\Lambda)$ is the unital $C^*$-algebra generated by $\{u_j\}_{j\in\Lambda}$.
That contradicts \cite[Proposition~4.1.1]{Brown-06}, since $\mathbb F_\Lambda$ is nonamenable.

By \Cref{lem compact-perturbation}, we have $\mathcal N=\mathcal P\oplus\mathcal Q\oplus\mathcal D$ with $\mathcal Q$ and $\mathcal D$ as in the statement.
By \Cref{lem nonzero-intertwiner}, $\mathcal N$ has a nonzero type $\mathrm{II}_1$ central summand.
Since this summand lies in $\mathcal P$, we have $\mathcal P\ne0$.
The required isomorphism $\theta$ is given by \Cref{lem compact-perturbation}.
\end{proof}

\section{Trace-class perturbations}
\label{sec trace-class}

In this section, we assume that $\mathcal H$ is the full Fock space $\mathcal F(\ell^2(\Lambda))$ and we will prove \Cref{thm main}.

\subsection{Finite truncations}

As in \eqref{equ sum-K}, let $\{K_j\}_{j\in\Lambda}$ be a family of self-adjoint trace-class operators in $\mathcal B(\mathcal H)$ such that
\begin{equation*}
	\gamma_j=\|K_j\|_1<\infty\quad(j\in\Lambda),\qquad \gamma=\sum_{j\in\Lambda}\gamma_j<\infty.
\end{equation*}
Then at most countably many $K_j$ are nonzero.
The spectral decompositions of the nonzero operators $K_j$ give vectors $\{\xi_m\}_{m=1}^\infty$, signs $\sigma_m\in\{-1,1\}$, and labels $i(m)\in\Lambda$ such that
\begin{equation*}
    K_j=\sum_{i(m)=j}\sigma_m(\xi_m\,\hat\otimes\,\xi_m),
    \qquad
    \gamma_j=\sum_{i(m)=j}\|\xi_m\|_2^2,
\end{equation*}
where $\{\xi_m\colon i(m)=j\}$ is an orthogonal set for each $j\in\Lambda$.
If there are only finitely many terms in total, append zero vectors with arbitrary signs and indices.
Moreover, we have $\gamma=\sum_{m=1}^\infty\|\xi_m\|_2^2<\infty$.

Let $e_j=l_j\Omega=r_j\Omega$ for $j\in\Lambda$.
These vectors form the standard orthonormal basis of $\ell^2(\Lambda)$.
For $q\geqslant1$ and $(j_1,\ldots,j_q)\in\Lambda^q$, write
$w=j_1\cdots j_q$ and define the corresponding \emph{Fock word} by
\begin{equation*}
    e_w=e_{j_1}\otimes\cdots\otimes e_{j_q}.
\end{equation*}
For the empty word $\varnothing$, let $e_{\varnothing}=\Omega$.
The vectors $e_w$, indexed by all finite words in $\Lambda$ including the empty word, form an orthonormal basis of $\mathcal H$.

\begin{lemma}\label{lem polynomial-vector}
The linear span of the Fock words consists precisely of the vectors $p\Omega$, where $p$ is a noncommutative polynomial in $\{s_j\}_{j\in\Lambda}$.
This span is preserved by $r_j$ and $r_j^*$ for every $j\in\Lambda$.
\end{lemma}

\begin{proof}
Clearly, this span contains $\Omega$ and is invariant under each $s_j=l_j+l_j^*$, and hence contains every $p\Omega$.
Conversely, $e_{\varnothing}=\Omega$ and $e_j=s_j\Omega$ for $j\in\Lambda$.
For every $i,j\in\Lambda$ and every word $v$, possibly empty, we have $l_ie_{jv}=e_{ijv}$ and $l_i^*e_{jv}=\delta_{ij}e_v$.
Hence
\begin{equation*}
    e_{ijv}=s_i e_{jv}-\delta_{ij}e_v.
\end{equation*}
Induction on word length now gives $e_w=p_w\Omega$ for some polynomial $p_w$.
Preservation by right creation and annihilation follows by adding or removing the last letter of a Fock word.
\end{proof}

Fix $n\geqslant 1$.
For each $1\leqslant m\leqslant n$, choose a polynomial $a_{m,n}$ in $\{s_j\}_{j\in\Lambda}$ such that $a_{m,n}\Omega$ is a finite truncation of the Fock expansion of $\xi_m$ with $\|a_{m,n}\Omega-\xi_m\|\leqslant n^{-1}\|\xi_m\|$.
Let $a_{m,n}=0$ for $m>n$.
It is clear that $a_{m,n}\in\mathcal L(\mathbb F_\Lambda)$ and
\begin{equation}\label{equ vector-approx}
	\|a_{m,n}\Omega\|\leqslant\|\xi_m\|,
	\qquad
	\lim_{n\to\infty}\sum_{m=1}^\infty\|a_{m,n}\Omega-\xi_m\|^2=0.
\end{equation}
Define
\begin{equation}\label{equ K-jn}
    K_{j,n}=\sum_{\substack{1\leqslant m\leqslant n\\ i(m)=j}}
    \sigma_m(a_{m,n}\Omega\,\hat\otimes\,a_{m,n}\Omega),
    \qquad
    S_{j,n}=s_j+K_{j,n}.
\end{equation}
Then $\lim_{n\to\infty}\sum_{j\in\Lambda}\|K_{j,n}-K_j\|_1=0$.
For every $n\geqslant1$, we have $K_{j,n}=0$ except for finitely many $j\in\Lambda$.
More precisely, $K_{j,n}=0$ whenever $j\notin\{i(1),i(2),\ldots,i(n)\}$.

By \Cref{prop compact-perturbation}, to obtain the decomposition in \Cref{thm main}, it suffices to construct a nonzero $Y\in\mathcal B(\mathcal H)$ such that $S_jY=Ys_j$ for every $j\in\Lambda$.
We first study the intertwining equations for the finite-rank perturbations $S_{j,n}$.

\subsection{Intertwiners for finite-rank perturbations}

Fix $n\geqslant1$.
For $1\leqslant m\leqslant n$, let $a_{m,n}\in\mathcal L(\mathbb F_\Lambda)$ be the finite truncation of the Fock expansion of $\xi_m$ chosen in the preceding subsection.
By \eqref{equ K-jn}, the equation $S_{j,n}Y=Ys_j$ is equivalent to
\begin{equation*}
	[s_j,Y]
	=-\sum_{\substack{1\leqslant m\leqslant n\\i(m)=j}}
	\sigma_m(a_{m,n}\Omega\,\hat{\otimes}\,Y^*a_{m,n}\Omega).
\end{equation*}
For $A_m\in\mathcal L(\mathbb F_\Lambda)'$, the identity $[s_j,r_i^*]=-\delta_{ij}P_\Omega$ gives
\begin{equation*}
	[s_j,R_{a_{m,n}}r_{i(m)}^*A_m]
	=R_{a_{m,n}}[s_j,r_{i(m)}^*]A_m
	=-\delta_{i(m),j}(a_{m,n}\Omega\,\hat{\otimes}\,A_m^*\Omega).
\end{equation*}
This suggests the candidate expression
\begin{equation}\label{equ Y=WA}
	Y=A_0+\sum_{m=1}^nR_{a_{m,n}}r_{i(m)}^*A_m,
	\qquad
	A_0,A_1,\ldots,A_n\in \mathcal L(\mathbb F_\Lambda)'.
\end{equation}
Let
\begin{equation*}
	\mathscr H=\bigoplus_{m=0}^{\infty}\mathcal H,
	\qquad
	\mathfrak A=(I\otimes \mathcal L(\mathbb F_\Lambda))'=\mathcal B(\ell^2)\,\overline\otimes\,\mathcal L(\mathbb F_\Lambda)',
	\qquad
	\mathcal T'=\operatorname{Tr}\otimes\tau'.
\end{equation*}
Let $J_m\colon\mathcal H\to\mathscr H$ be the $m$-th coordinate embedding and let $\Omega_m=J_m\Omega$ for $m\geqslant 0$.
Write $E_n=\sum_{m=0}^nJ_mJ_m^*$.
Define
\begin{equation}\label{equ Wn}
    W_n=\bigl(I,R_{a_{1,n}}r_{i(1)}^*,\ldots,R_{a_{n,n}}r_{i(n)}^*\bigr)
    \colon E_n\mathscr H\to\mathcal H
\end{equation}
and extend $W_n$ by zero on $(1-E_n)\mathscr H$.
Define $\widetilde A_n\colon\mathcal H\to E_n\mathscr H$ by $\widetilde A_n=\sum_{m=0}^n J_mA_m$.
Then \eqref{equ Y=WA} becomes $Y=W_n\widetilde A_n$.
Write $\widehat s_j=I\otimes s_j$ in the following lemma.

\begin{lemma}\label{lem truncated}
For every $n\geqslant 1$, there exists a projection $P_n\in E_n\mathfrak A E_n$ satisfying
\begin{align}
	\mathcal T'(P_n)&\geqslant1,\label{equ ker-dim}\\
	(S_{j,n}W_n-W_n\widehat s_j)P_n&=0
	\qquad(j\in\Lambda),\label{equ truncated-1}\\
	P_n(W_n^*a_{m,n}\Omega-\sigma_m\Omega_m) &=0
	\qquad(m\geqslant1).\label{equ truncated-2}
\end{align}
\end{lemma}

\begin{proof}
Let $\eta_m=W_n^*a_{m,n}\Omega-\sigma_m\Omega_m$ for $1\leqslant m\leqslant n$, let $\mathscr K_n$ be the closed submodule of $E_n\mathscr H$ generated by $\eta_1,\eta_2,\ldots,\eta_n$ under the diagonal left action of
$\mathcal L(\mathbb F_\Lambda)$, and let $P_n$ be the projection onto
$E_n\mathscr H\ominus\mathscr K_n$.
Then $P_n\in E_n\mathfrak A E_n$, and \Cref{lem cyclic-dimension} gives
\begin{equation*}
    \mathcal T'(P_n)
    =n+1-\dim_{\mathcal L(\mathbb F_\Lambda)}\mathscr K_n
    \geqslant1.
\end{equation*}
Clearly, $P_n\eta_m=0$ for $1\leqslant m\leqslant n$.
For $m>n$, we have $a_{m,n}=0$ and $P_n\Omega_m=0$.
This proves \eqref{equ truncated-2}.
For $y\in P_n\mathscr H$, \eqref{equ K-jn} and the identity $[s_j,r_i^*]=-\delta_{ij}P_\Omega$ give
\begin{equation*}
    (S_{j,n}W_n-W_n\widehat s_j)y
    =\sum_{\substack{1\leqslant m\leqslant n\\i(m)=j}}
      \sigma_m\langle y,\eta_m\rangle\,a_{m,n}\Omega=0.
\end{equation*}
This proves \eqref{equ truncated-1}.
\end{proof}

To control the passage to a weak limit, we use $\varphi$ to estimate $\mathcal T'(Y^*Y)$ for rows of the form $Y=W_nB$.

\begin{lemma}\label{lem rowtrace}
For $1\leqslant m\leqslant n$, let $a_m\in\mathcal L(\mathbb F_\Lambda)$ and $i(m)\in\Lambda$.
Let
\begin{equation*}
    W_n=\bigl(I,R_{a_1}r_{i(1)}^*,\ldots,R_{a_n}r_{i(n)}^*\bigr),
    \qquad c=\left(\sum_{m=1}^n\|a_m\Omega\|^2\right)^{1/2}.
\end{equation*}
Suppose that $B\in E_n\mathfrak A E_n$, $Y=W_nB$, and $Y^*Y\in E_n\mathfrak A E_n$. Then $YY^*\in\mathcal E$ and
\begin{equation*}
    \mathcal T'(Y^*Y)=\varphi(YY^*)
    \leqslant\|Y^*\Omega\|^2+
        2c\|Y\|\,\mathcal T'(B^*B)^{1/2}.
\end{equation*}
\end{lemma}

\begin{proof}
Write $Y=(C_0,C_1,\ldots,C_n)$ on $E_n\mathscr H$, where $C_k=YJ_k$.
Each entry of $W_n$ belongs to $\mathcal E$, and each entry of $B$ belongs to $\mathcal L(\mathbb F_\Lambda)'$.
Since $Y=W_nB$, we have $C_k\in\mathcal E$ for every $0\leqslant k\leqslant n$.
Thus, $YY^*=\sum_{k=0}^n C_kC_k^*\in\mathcal E$.

Since $Y^*Y\in E_n\mathfrak A E_n$, each diagonal entry
$C_k^*C_k=J_k^*Y^*YJ_k$ belongs to $\mathcal L(\mathbb F_\Lambda)'$.
By \Cref{prop full-trace}, we have $\varphi(C_kC_k^*)=\varphi(C_k^*C_k)=\tau'(C_k^*C_k)$.
Summing over $k$ gives
\begin{equation*}
    \varphi(YY^*)=\sum_{k=0}^n\tau'(C_k^*C_k)=\mathcal T'(Y^*Y).
\end{equation*}
Let $\partial_jY=s_jY-Y\widehat s_j$ for $j\in\Lambda$.
Since $B$ commutes with $\widehat s_j=I\otimes s_j$, it follows that $\partial_jY=(s_jW_n-W_n\widehat s_j)B$.
Using $[s_j,r_i^*]=-\delta_{ij}P_\Omega$ and $R_{a_m}\Omega=a_m\Omega$, we obtain
\begin{equation*}
	\partial_jY=-\sum_{\substack{1\leqslant m\leqslant n\\i(m)=j}}
	a_m\Omega\,\hat{\otimes}\,B^*\Omega_m.
\end{equation*}
The Cauchy--Schwarz inequality therefore gives
\begin{equation*}
	\sum_{j\in\Lambda}\|\partial_jY\|_1
	\leqslant
	\sum_{m=1}^n\|a_m\Omega\|\|B^*\Omega_m\|\leqslant c\left(\sum_{m=1}^n\|B^*\Omega_m\|^2\right)^{1/2}
	\leqslant
	c\,\mathcal T'(B^*B)^{1/2},
\end{equation*}
where the last inequality uses $\sum_{m=0}^n\|B^*\Omega_m\|^2=\mathcal T'(BB^*)=\mathcal T'(B^*B)$.
Since
\begin{equation*}
	[s_j,YY^*]=(\partial_jY)Y^*-Y(\partial_jY)^*,
\end{equation*}
we have $\|[s_j,YY^*]\|_1\leqslant 2\|Y\|\|\partial_jY\|_1$.
Since $\|l_j\|=1$, the definition of $\varphi$ gives
\begin{equation*}
    \bigl|\varphi(YY^*)-\|Y^*\Omega\|^2\bigr|
    \leqslant\sum_{j\in\Lambda}\|[s_j,YY^*]\|_1
    \leqslant2c\|Y\|\,\mathcal T'(B^*B)^{1/2}.
\end{equation*}
This proves the desired estimate.
\end{proof}

\subsection{Graph normalization and a nonzero weak limit}
\label{subsec graph}

Let $P_n$ be the projection in \Cref{lem truncated} and let $T_n=W_nP_n$.
Since $P_n$ commutes with $\widehat s_j$, it follows from \eqref{equ truncated-1} that
\begin{equation*}
	S_{j,n}T_n=T_n\widehat s_j\qquad(j\in\Lambda).
\end{equation*}
Taking adjoints gives $T_n^*S_{j,n}=\widehat s_jT_n^*$.
Thus, $T_n^*T_n$ commutes with every $\widehat s_j$, and hence $T_n^*T_n\in P_n\mathfrak AP_n$.
For the operator $T_n$, its \emph{graph map} $\begin{pmatrix}P_n\\T_n\end{pmatrix}\colon\mathscr H\to\mathscr H\oplus\mathcal H$ has polar decomposition
\begin{equation*}
	\begin{pmatrix}P_n\\T_n\end{pmatrix}
	=\begin{pmatrix}B_n\\Y_n\end{pmatrix}
	(P_n+T_n^*T_n)^{1/2},
\end{equation*}
where
\begin{equation}\label{equ B-Y-1}
	B_n=P_n(I+T_n^*T_n)^{-1/2}\in E_n\mathfrak AE_n,
	\qquad
	Y_n=W_nB_n,
	\qquad
	B_n^2+Y_n^*Y_n=P_n.
\end{equation}
In particular, $Y_n$ is a contraction and $Y_n^*Y_n\in E_n\mathfrak AE_n$.
Since $B_n$ commutes with every $\widehat s_j$, we also have
\begin{equation}\label{equ SY=Ys}
    S_{j,n}Y_n=Y_n\widehat s_j\qquad(j\in\Lambda).
\end{equation}
Since $B_n=B_nP_n$ and $Y_n^*=B_nW_n^*$, multiplying \eqref{equ truncated-2} on the left by $B_n$ gives
\begin{equation}\label{equ B-Y-2}
    Y_n^*a_{m,n}\Omega=\sigma_m B_n\Omega_m
    \qquad(m\geqslant1).
\end{equation}
Thus, $\|B_n\Omega_m\|\leqslant\|a_{m,n}\Omega\|\leqslant\|\xi_m\|$ by \eqref{equ vector-approx}.
Together with $\|B_n\Omega_0\|\leqslant1$, the vacuum trace on the diagonal entries of $B_n^2$ gives
\begin{equation}\label{equ g_n}
    g_n:=\mathcal T'(B_n^2)
    =\sum_{m=0}^{\infty}\|B_n\Omega_m\|^2
    \leqslant1+\gamma.
\end{equation}
We now pass to a weak-operator limit and show that the resulting intertwiner has initial support of trace at least $1$.

\begin{lemma}\label{lem intertwiner}
Let $Y_\alpha$ be the weak-operator limit of $\{Y_n\}_{n=1}^\infty$ along a free ultrafilter $\alpha$ on $\mathbb N$.
Then
\begin{equation*}
    S_jY_\alpha=Y_\alpha\widehat s_j\qquad(j\in\Lambda),
    \qquad
    \mathcal T'(P_\alpha)\geqslant1,
\end{equation*}
where $P_\alpha$ is the support of $Y_\alpha^*Y_\alpha\in\mathfrak A$.
\end{lemma}

\begin{proof}
Since $B_n$ and $Y_n$ are contractions, weak-operator compactness gives the limits
\begin{equation*}
    B_\alpha=\lim_{n\to\alpha}B_n\in\mathfrak A,
    \qquad
    Y_\alpha=\lim_{n\to\alpha}Y_n\in\mathcal B(\mathscr H,\mathcal H).
\end{equation*}
Since $S_{j,n}\to S_j$ in norm and $S_{j,n}Y_n=Y_n\widehat s_j$ by \eqref{equ SY=Ys}, we have
\begin{equation*}
    S_jY_\alpha=Y_\alpha\widehat s_j\qquad(j\in\Lambda).
\end{equation*}
Taking adjoints shows that $Y_\alpha^*Y_\alpha$ commutes with every $\widehat s_j$.
Thus, $Y_\alpha^*Y_\alpha\in\mathfrak A$, and hence its support $P_\alpha$ belongs to $\mathfrak A$.
Let $Q_\alpha=I-P_\alpha$.

Let $w_n=W_n^*\Omega=\bigl(\Omega,(r_{i(m)}a_{m,n}^*\Omega)_{m\geqslant1}\bigr)$ and $w=\bigl(\Omega,(r_{i(m)}\xi_m^*)_{m\geqslant1}\bigr)$ be vectors in $\mathscr H$, where $\xi\mapsto\xi^*$ is the conjugate-linear isometry on $\mathcal H$ extending $a\Omega\mapsto a^*\Omega$.
For every $m\geqslant 1$, we have $a_{m,n}\Omega\to\xi_m$ in norm.
Together with \eqref{equ vector-approx}, it follows that $w_n\to w$ in norm.
Since $Y_n^*a_{m,n}\Omega=\sigma_mB_n\Omega_m$ by \eqref{equ B-Y-2} and $Y_n^*\Omega=B_nw_n$, we have
\begin{equation*}
    B_\alpha\Omega_m=\sigma_mY_\alpha^*\xi_m\quad(m\geqslant1),
    \qquad
    B_\alpha w=Y_\alpha^*\Omega.
\end{equation*}
Then $Q_\alpha B_\alpha\Omega_m=0$ for $m\geqslant1$ and $Q_\alpha B_\alpha w=0$.
Since $Q_\alpha B_\alpha\in\mathfrak A$ commutes with the diagonal left action of $\mathcal L(\mathbb F_\Lambda)$ and $\mathcal L(\mathbb F_\Lambda)\Omega$ is dense in $\mathcal H$, it follows that $Q_\alpha B_\alpha=0$.
Since $0\leqslant B_n\leqslant P_n$, we have
\begin{equation*}
    \lim_{n\to\alpha}\|B_nQ_\alpha\eta\|^2
    \leqslant\langle B_\alpha Q_\alpha\eta,Q_\alpha\eta\rangle=0
    \qquad(\eta\in\mathscr H).
\end{equation*}
Hence $B_nQ_\alpha\to 0$ in the strong-operator topology along $\alpha$.

Suppose to the contrary that $\mathcal T'(P_\alpha)=\sum_{m=0}^\infty\|P_\alpha\Omega_m\|^2<1$ and let
\begin{equation*}
    h_n=\mathcal T'(Q_\alpha B_n^2Q_\alpha)
    =\sum_{m=0}^\infty\|B_nQ_\alpha\Omega_m\|^2.
\end{equation*}
Then each summand tends to zero along $\alpha$.
Since $\|B_n\Omega_m\|=\|Y_n^*a_{m,n}\Omega\|\leqslant\|\xi_m\|$ by \eqref{equ vector-approx} and \eqref{equ B-Y-2}, we have
\begin{equation*}
    \|B_nQ_\alpha\Omega_m\|
    \leqslant\|B_n\Omega_m\|+\|B_nP_\alpha\Omega_m\|
    \leqslant\|\xi_m\|+\|P_\alpha\Omega_m\|
    \qquad(m\geqslant1).
\end{equation*}
Recall that $\gamma=\sum_{m=1}^\infty\|\xi_m\|^2<\infty$.
Thus, $\lim_{n\to\alpha}h_n=0$ by Lebesgue's dominated convergence theorem.

Since $h_n=\mathcal T'(B_nQ_\alpha B_n)=\sum_{m=0}^\infty\|Q_\alpha B_n\Omega_m\|^2$, we have $Q_\alpha B_n\Omega_m\to0$ along $\alpha$ for every $m\geqslant0$.
Since every $Q_\alpha B_n\in\mathfrak A$ is a contraction commuting with the diagonal left action of $\mathcal L(\mathbb F_\Lambda)$ and $\mathcal L(\mathbb F_\Lambda)\Omega$ is dense in $\mathcal H$, it follows that $Q_\alpha B_n\to 0$ in the strong-operator topology along $\alpha$.
Combining this with $w_n\to w$ as $n\to\infty$, we obtain $\lim_{n\to\alpha}\|Q_\alpha B_nw_n\|=0$.

Let $N\geqslant n\geqslant 1$.
By $Y_n^*\Omega=B_nw_n$, $\sum_{m=1}^n\|a_{m,n}\Omega\|^2\leqslant\gamma$, \eqref{equ ker-dim}, \eqref{equ B-Y-1}, and applying \Cref{lem rowtrace} to $Y_nQ_\alpha E_N=W_n(B_nQ_\alpha E_N)$ with $E_N=\sum_{m=0}^N J_mJ_m^*$, we obtain
\begin{equation*}
    \mathcal T'(E_NQ_\alpha P_nQ_\alpha E_N)
    \leqslant h_n+\|Q_\alpha B_nw_n\|^2+2\sqrt{\gamma h_n}.
\end{equation*}
Here we used $\|Y_nQ_\alpha E_N\|\leqslant1$ and $\mathcal T'(E_NQ_\alpha B_n^2Q_\alpha E_N)\leqslant h_n$.
By letting $N\to\infty$, it follows that
\begin{equation*}
	\mathcal T'(P_nQ_\alpha)=\mathcal T'(Q_\alpha P_nQ_\alpha)
    \leqslant h_n+\|Q_\alpha B_nw_n\|^2+2\sqrt{\gamma h_n}.
\end{equation*}
Consequently, $\lim_{n\to\alpha}\mathcal T'(P_nQ_\alpha)=0$.
On the other hand, by $Q_\alpha=I-P_\alpha$ and \eqref{equ ker-dim}, we have
\begin{equation*}
    \mathcal T'(P_nQ_\alpha)
    =\mathcal T'(P_n)-\mathcal T'(P_\alpha P_nP_\alpha)
    \geqslant1-\mathcal T'(P_\alpha)>0.
\end{equation*}
That is a contradiction.
Therefore, $\mathcal T'(P_\alpha)\geqslant1$.
\end{proof}

\begin{proof}[Proof of \Cref{thm main}]
By \Cref{lem intertwiner}, the operator $Y_\alpha$ is nonzero and satisfies
\begin{equation*}
    S_jY_\alpha=Y_\alpha\widehat s_j\qquad(j\in\Lambda).
\end{equation*}
Choose $m\geqslant0$ such that $Y=Y_\alpha J_m\neq0$.
Using $\widehat s_jJ_m=J_ms_j$, we obtain
\begin{equation*}
	S_jY=S_jY_\alpha J_m=Y_\alpha\widehat s_jJ_m=Ys_j
	\qquad(j\in\Lambda).
\end{equation*}
The required decomposition and the isomorphism follow from \Cref{prop compact-perturbation}.

Suppose that $\Lambda$ is infinite and choose distinct indices $\{j_n\}_{n=1}^\infty$ in $\Lambda$.
For every unit vector $\xi\in\mathcal H$, we have $\sum_{n=1}^\infty\|l_{j_n}^*\xi\|^2\leqslant1$ and $\sum_{n=1}^\infty\|K_{j_n}\|_1\leqslant\gamma$.
It follows that
\begin{equation*}
    \|S_{j_n}\xi-l_{j_n}\xi\|
    \leqslant\|l_{j_n}^*\xi\|+\|K_{j_n}\|_1\to0,
\end{equation*}
Since $\{l_{j_n}\xi\}_{n=1}^\infty$ is an orthonormal set, $\{S_{j_n}\xi\}_{n=1}^\infty$ has no convergent subsequence.
Thus, no nonzero finite-dimensional subspace of $\mathcal H$ is invariant under every $S_j$.
As a consequence, we have $\mathcal D=0$.
\end{proof}

Taking two semicircular generators as the real and imaginary parts of a single operator gives a counterexample to the trace-class density of irreducible operators.

\begin{corollary}\label{cor irreducible}
Let $\Lambda=\{1,2\}$ and $T=s_1+\mathbf is_2$.
Then $T+K$ is reducible for every trace-class operator $K\in\mathcal B(\mathcal H)$. 
In particular, the set of irreducible operators in $\mathcal B(\mathcal H)$ is not dense in the trace-class norm.
\end{corollary}

\begin{proof}
Write $K=K_1+\mathbf iK_2$, where $K_1$ and $K_2$ are self-adjoint trace-class operators.
Then
\begin{equation*}
	W^*(T+K)=W^*(s_1+K_1,s_2+K_2).
\end{equation*}
Since the type $\mathrm{II}_1$ summand of $W^*(T+K)$ is isomorphic to $\mathcal L(\mathbb F_2)$ by \Cref{thm main}, $T+K$ is not an irreducible operator in $\mathcal B(\mathcal H)$.
\end{proof}

\section{Murray-von Neumann dimension}
\label{sec multiplicity}

In this section, we study the Murray-von Neumann dimension and prove \Cref{thm multiplicity-range}.
For a self-adjoint trace-class family $\{K_j\}_{j\in\Lambda}$ as in \Cref{thm main}, let
\begin{equation*}
	d(\{K_j\}_{j\in\Lambda})=\dim_{\mathcal P}(z_{\mathcal P}\mathcal H).
\end{equation*}
Since $\mathcal P$ is isomorphic to $\mathcal L(\mathbb F_\Lambda)$, the dimension can be viewed as the multiplicity of $\mathcal L(\mathbb F_\Lambda)$ under the perturbation by $\{K_j\}_{j\in\Lambda}$.
We bound this dimension from below using \Cref{lem intertwiner} and from above by estimating $\varphi(z_{\mathcal P})$.

\begin{lemma}\label{lem multiplicity}
Under the hypotheses of \Cref{thm main}, if $\gamma>0$, then
\begin{equation*}
	1\leqslant d(\{K_j\}_{j\in\Lambda})<1+\gamma.
\end{equation*}
\end{lemma}

\begin{proof}
Let $Y_\alpha$ be the intertwiner in \Cref{lem intertwiner} and let $Y_\alpha=V_\alpha|Y_\alpha|$ be its polar decomposition.
Then $\mathcal T'(P_\alpha)\geqslant 1$ and $P_\alpha=V_\alpha^*V_\alpha$.
As in the proof of \Cref{lem nonzero-intertwiner}, $V_\alpha$ identifies $P_\alpha\mathscr H$ with $V_\alpha V_\alpha^*\mathcal H$, and the compression of $\mathcal N$ to $V_\alpha V_\alpha^*\mathcal H$ is isomorphic to $\mathcal L(\mathbb F_\Lambda)$.
It follows from \Cref{thm main} that $V_\alpha V_\alpha^*\leqslant z_{\mathcal P}$, and hence
\begin{equation*}
    d(\{K_j\}_{j\in\Lambda})\geqslant\dim_{\mathcal P}(V_\alpha V_\alpha^*\mathcal H)
    =\mathcal T'(P_\alpha)\geqslant1.
\end{equation*}
We use the notation in \Cref{subsec tracial}.
For $X\in\mathcal N'$, we have $[s_j,X]=[X,K_j]$ for $j\in\Lambda$.
It follows that $X\in\mathcal E$ and
\begin{equation*}
    \varphi(X)=\langle X\Omega,\Omega\rangle+\operatorname{Tr}(DX),
    \qquad
    D=\sum_{j\in\Lambda}[K_j,l_j]\in\mathcal S_1(\mathcal H).
\end{equation*}
Thus, $\varphi$ is normal on $\mathcal N'$.
Suppose that $Y\in\mathcal B(\mathcal H)$ satisfies $S_jY=Ys_j$ for every $j\in\Lambda$.
Then $Y^*Y\in\mathcal L(\mathbb F_\Lambda)'$ and $[s_j,Y]=-K_jY$.
In particular, $Y\in\mathcal E$.
By \Cref{prop full-trace}, we have
\begin{equation}\label{equ phi-Y}
    \varphi(YY^*)=\varphi(Y^*Y)=\tau'(Y^*Y).
\end{equation}
Via the isomorphism $\theta$ in \Cref{thm main}, we can decompose $z_{\mathcal P}\mathcal H$ as an orthogonal sum of cyclic $\mathcal L(\mathbb F_\Lambda)$-modules.
Thus, there are partial isometries $\{V_k\}$ in $\mathcal B(\mathcal H)$ such that
\begin{equation*}
	S_jV_k=V_ks_j\quad(j\in\Lambda),
	\qquad
	z_{\mathcal P}=\sum_kV_kV_k^*.
\end{equation*}
By \eqref{equ phi-Y} and the normality of $\varphi$ on $\mathcal N'$, we have
\begin{equation*}
    d(\{K_j\}_{j\in\Lambda})=\sum_k\tau'(V_k^*V_k)=\sum_k\varphi(V_kV_k^*)=\varphi(z_{\mathcal P}).
\end{equation*}
For $j\in\Lambda$, let $C_j=[s_j,z_{\mathcal P}]=[z_{\mathcal P},K_j]$.
Then $C_j$ is skew-adjoint and
\begin{equation*}
    \|C_j\|_1=\frac12\|[2z_{\mathcal P}-1,K_j]\|_1\leqslant\|K_j\|_1.
\end{equation*}
Since $|\langle l_j\xi,\xi\rangle|<\|\xi\|^2$ for every nonzero $\xi\in\mathcal H$, we have $|\operatorname{Tr}(C_jl_j)|<\|C_j\|_1$ whenever $C_j\neq0$.
It follows that $|\operatorname{Tr}(C_jl_j)|<\|K_j\|_1$ whenever $K_j\neq0$.
Since $\gamma>0$, at least one $K_j$ is nonzero. Hence
\begin{equation*}
    d(\{K_j\}_{j\in\Lambda})=\varphi(z_{\mathcal P})
    \leqslant\|z_{\mathcal P}\Omega\|^2+\sum_{j\in\Lambda}|\operatorname{Tr}(C_jl_j)|
    <1+\gamma.
\end{equation*}
This completes the proof.
\end{proof}

Choose an index in $\Lambda$ and denote it by $1$.
For any self-adjoint trace-class operator $L$ on $\mathcal H$, write
\begin{equation*}
	d(L)=d(\{K_j\}_{j\in\Lambda}),
\end{equation*}
where $K_1=L$ and $K_j=0$ for every $j\in\Lambda\setminus\{1\}$.
For $a\in\mathcal L(\mathbb F_\Lambda)$ and $\lambda>0$, let
\begin{equation}\label{equ rank-one-data}
    b=(\lambda-s_1)a,
    \qquad
    L=\lambda^{-1}b\Omega\,\hat{\otimes}\,b\Omega,
    \qquad
    S_j=s_j+\delta_{j1}L\quad(j\in\Lambda).
\end{equation}
Let $\mathcal N=W^*(\{S_j\colon j\in\Lambda\})$, and let $\mathcal P=z_{\mathcal P}\mathcal N$ be the summand given by \Cref{thm main}.
For the rank-one perturbations in \eqref{equ rank-one-data} we establish an upper bound for $d(L)$.

\begin{lemma}\label{lem rank-one-upper}
Let $p\in\mathcal L(\mathbb F_\Lambda)$ be a nonzero projection.
With the notation of \eqref{equ rank-one-data}, suppose that $a$ is an invertible operator in $p\mathcal L(\mathbb F_\Lambda)p$ such that $\langle a\Omega,b\Omega\rangle=\lambda$.
Then $R_p\in\mathcal N'$ and
\begin{equation*}
    d(L)=\tau(1-p)+\dim_{\mathcal P}(z_{\mathcal P}R_p\mathcal H)
    \leqslant1+\tau(1-p).
\end{equation*}
\end{lemma}

\begin{proof}
Since $b=bp$, we have $b\Omega\in R_p\mathcal H$ and $R_pL=LR_p=L$.
Thus, $R_p$ commutes with every $S_j$, and hence $R_p\in\mathcal N'$.
On $R_{1-p}\mathcal H$, the family $\{S_j\}_{j\in\Lambda}$ agrees with $\{s_j\}_{j\in\Lambda}$.
By \Cref{thm main}, we have $R_{1-p}\mathcal H\subseteq z_{\mathcal P}\mathcal H$ and
\begin{equation*}
    \dim_{\mathcal P}(R_{1-p}\mathcal H)
    =\dim_{\mathcal L(\mathbb F_\Lambda)}(R_{1-p}\mathcal H)
    =\tau(1-p).
\end{equation*}
Since $R_p\in\mathcal N'$ and $R_{1-p}\mathcal H\subseteq z_{\mathcal P}\mathcal H$, we have $z_{\mathcal P}\mathcal H=R_{1-p}\mathcal H\oplus z_{\mathcal P}R_p\mathcal H$.
It remains to prove that $\dim_{\mathcal P}(z_{\mathcal P}R_p\mathcal H)\leqslant1$.

Let $\mathcal K=\overline{\mathcal Na\Omega+\mathcal Nb\Omega}$.
Since $a\Omega,b\Omega\in R_p\mathcal H$ and $R_p\in\mathcal N'$, we have $\mathcal K\subseteq R_p\mathcal H$.
Since $\mathcal K$ is invariant under every $S_j$ and contains $\operatorname{ran}L$, it is invariant under every $s_j=S_j-\delta_{j1}L$.
By the invertibility of $a$ in $p\mathcal L(\mathbb F_\Lambda)p$, we have
\begin{equation*}
    R_p\mathcal H
    =\overline{\mathcal L(\mathbb F_\Lambda)p\Omega}
    =\overline{\mathcal L(\mathbb F_\Lambda)a\Omega}
    \subseteq\mathcal K.
\end{equation*}
Thus, $\mathcal K=R_p\mathcal H$.
Since $\langle a\Omega,b\Omega\rangle=\lambda$, we have $La\Omega=b\Omega$, and hence
\begin{equation}\label{equ S1-eigen}
    S_1a\Omega=s_1a\Omega+b\Omega=\lambda a\Omega.
\end{equation}
Under the isomorphism $\theta$ in \Cref{thm main}, $z_{\mathcal P}S_1$ corresponds to $s_1$, which has no eigenvalues.
It follows from \eqref{equ S1-eigen} that $z_{\mathcal P}a\Omega=0$, and therefore
\begin{equation*}
    z_{\mathcal P}R_p\mathcal H
    =z_{\mathcal P}\mathcal K
    =\overline{\mathcal Pz_{\mathcal P}b\Omega}.
\end{equation*}
By \Cref{lem cyclic-dimension}, this cyclic $\mathcal P$-module has dimension at most $1$.
\end{proof}

When the projection $p$ in \Cref{lem rank-one-upper} commutes with $r_1$, we can compute $d(L)$ explicitly.

\begin{lemma}\label{lem rank-one-dimension}
Under the hypotheses of \Cref{lem rank-one-upper}, suppose that
$p$ commutes with $r_1$. Then
\begin{equation*}
    \dim_{\mathcal P}(z_{\mathcal P}R_p\mathcal H)=1,
    \qquad d(L)=1+\tau(1-p).
\end{equation*}
\end{lemma}

\begin{proof}
For $x\in\mathcal L(\mathbb F_\Lambda)$, let
$\widetilde x=(x|_{R_p\mathcal H})\oplus x$ be the diagonal action on
$R_p\mathcal H\oplus\mathcal H$. Define
\begin{equation*}
    W\colon R_p\mathcal H\oplus\mathcal H\to R_p\mathcal H,
    \qquad W(y_0,y_1)=y_0+R_b r_1y_1.
\end{equation*}
Let $\xi=W^*a\Omega\in R_p\mathcal H\oplus\mathcal H$ and let
\begin{equation*}
    \mathscr L
    =\{\widetilde x\xi\colon x\in\mathcal L(\mathbb F_\Lambda)\}^\perp
    \subseteq R_p\mathcal H\oplus\mathcal H.
\end{equation*}
Since $p$ commutes with $r_1$, we have $pW=W\widetilde p$.
Together with $pa=a$, this gives $\widetilde p\xi=\xi$.
Since $\mathscr L^\perp$ is generated by $\xi$ under the diagonal action
and $\widetilde p\xi=\xi$, \Cref{lem cyclic-dimension} gives
$\dim_{\mathcal L(\mathbb F_\Lambda)}\mathscr L^\perp\leqslant\tau(p)$.
It follows that
\begin{equation*}
    \dim_{\mathcal L(\mathbb F_\Lambda)}\mathscr L
    =1+\tau(p)-\dim_{\mathcal L(\mathbb F_\Lambda)}\mathscr L^\perp
    \geqslant1.
\end{equation*}

For $y=(y_0,y_1)\in R_p\mathcal H\oplus\mathcal H$, the identity $[s_j,r_1]=\delta_{j1}P_\Omega$ gives
\begin{equation*}
    (s_jW-W\widetilde s_j)y
    =\delta_{j1}\langle y_1,\Omega\rangle b\Omega.
\end{equation*}
Suppose that $y\in\mathscr L$.
Since $\mathscr L$ is invariant under the diagonal action, both $y$ and $\widetilde s_1y$ are orthogonal to $\xi=W^*a\Omega$; equivalently, both $Wy$ and $W\widetilde s_1y$ are orthogonal to $a\Omega$.
Using $b=(\lambda-s_1)a$, we obtain
\begin{equation}\label{equ Wy-b}
    \langle Wy,b\Omega\rangle
    =-\langle(s_1W-W\widetilde s_1)y,a\Omega\rangle
    =-\lambda\langle y_1,\Omega\rangle.
\end{equation}
It follows that $LWy=-\langle y_1,\Omega\rangle b\Omega$.
Since $S_j=s_j+\delta_{j1}L$ and $\mathscr L$ is invariant under every $\widetilde s_j$, we obtain
\begin{equation}\label{equ SW=Ws}
    S_j(W|_{\mathscr L})
    =(W|_{\mathscr L})(\widetilde s_j|_{\mathscr L})
    \qquad(j\in\Lambda).
\end{equation}
Suppose that $y\in\ker(W|_{\mathscr L})$.
By \eqref{equ SW=Ws}, this kernel is invariant under every $\widetilde s_j$, and hence under the diagonal action of $\mathcal L(\mathbb F_\Lambda)$. Applying \eqref{equ Wy-b} to $\widetilde x y$, we obtain
\begin{equation*}
    0=\langle W\widetilde x y,b\Omega\rangle
    =-\lambda\langle y_1,x^*\Omega\rangle
    \qquad(x\in\mathcal L(\mathbb F_\Lambda)).
\end{equation*}
Since $\mathcal L(\mathbb F_\Lambda)\Omega$ is dense in $\mathcal H$, it follows that $y_1=0$, and hence $y_0=Wy=0$.
Thus, $W|_{\mathscr L}$ is injective.
Combining this with \eqref{equ SW=Ws}, the polar decomposition of $W|_{\mathscr L}$ gives a unitary operator $U\colon\mathscr L\to\overline{W\mathscr L}$ satisfying
\begin{equation*}
	S_jU=U(\widetilde s_j|_{\mathscr L})
    \qquad(j\in\Lambda).
\end{equation*}
By \Cref{thm main}, $\overline{W\mathscr L}\subseteq z_{\mathcal P}R_p\mathcal H$, and hence
\begin{equation*}
    \dim_{\mathcal P}(z_{\mathcal P}R_p\mathcal H)
    \geqslant\dim_{\mathcal P}\overline{W\mathscr L}
    =\dim_{\mathcal L(\mathbb F_\Lambda)}\mathscr L
    \geqslant1.
\end{equation*}
We complete the proof by applying \Cref{lem rank-one-upper}.
\end{proof}

We now use \Cref{lem rank-one-dimension} to construct positive rank-one perturbations for each prescribed increase $0\leqslant h<1$ in dimension.

\begin{lemma}\label{lem small-gain}
Let $0\leqslant h<1$ and $c>2(1-\sqrt{1-h})$.
Then there exists a positive rank-one operator $L$ on $\mathcal H$ such that
\begin{equation*}
    d(L)=1+h,\qquad \|L\|_1=c.
\end{equation*}
\end{lemma}

\begin{proof}
Since $|\Lambda|\geqslant 2$, we may choose an index in $\Lambda\setminus\{1\}$ and denote it by $2$.
Let $p$ be a projection in $W^*(s_2)$ with $\tau(p)=1-h$.
By Propositions~12.11 and~12.13 and Theorem~14.10 in \cite{NS-06-book}, the spectrum of $ps_1^2p$ in $p\mathcal L(\mathbb F_\Lambda)p$ is
\begin{equation*}
    \bigl[(1-\sqrt{1-h})^2,(1+\sqrt{1-h})^2\bigr].
\end{equation*}
Thus, $e=\chi_{(0,c^2/4)}(ps_1^2p)\in p\mathcal L(\mathbb F_\Lambda)p$ is a nonzero spectral projection of $ps_1^2p$ and we have
$\|s_1e\Omega\|^2=\tau(ps_1^2pe)<(c^2/4)\tau(e)$.
For sufficiently small $\varepsilon>0$, the positive invertible operator
\begin{equation*}
    a=\frac{e+\varepsilon p}{\|(e+\varepsilon p)\Omega\|}
    \in p\mathcal L(\mathbb F_\Lambda)p
\end{equation*}
satisfies $\|a\Omega\|=1$ and $\|s_1a\Omega\|<c/2$.
The trace-preserving automorphism of $\mathcal L(\mathbb F_\Lambda)$ sending $s_1$ to $-s_1$ and fixing the other generators fixes $p$ and $s_1^2$, and hence fixes $a$.
It follows that $\langle s_1a\Omega,a\Omega\rangle=\tau(as_1a)=0$.
Let
\begin{equation*}
    \lambda=\frac{c+\sqrt{c^2-4\|s_1a\Omega\|^2}}{2},
\end{equation*}
and define $b$ and $L$ by \eqref{equ rank-one-data}.
Then $\langle a\Omega,b\Omega\rangle=\lambda>0$.
Thus, $L$ is positive of rank one with
\begin{equation*}
	\|L\|_1=\lambda^{-1}\|b\Omega\|^2=\lambda\|a\Omega\|^2+\lambda^{-1}\|s_1a\Omega\|^2=c.
\end{equation*}
Since $r_1$ commutes with $s_2$, it commutes with $p$.
It follows from \Cref{lem rank-one-dimension} that
\begin{equation*}
	d(L)=1+\tau(1-p)=1+h.
\end{equation*}
This completes the proof.
\end{proof}

By iterating the rank-one construction for a single generator, we realize every dimension in $[1,1+\gamma)$ at total trace norm $\gamma$.

\begin{proposition}\label{prop realization}
Let $\gamma>0$ and $1\leqslant t<1+\gamma$.
Then there exists a positive finite-rank operator $K_1$ such that, with $K_j=0$ for $j\in\Lambda\setminus\{1\}$, we have
\begin{equation*}
	\gamma=\|K_1\|_1=\sum_{j\in\Lambda}\|K_j\|_1,\qquad d(\{K_j\}_{j\in\Lambda})=t.
\end{equation*}
\end{proposition}

\begin{proof}
Let $h=t-1\geqslant 0$.
Since $h<\gamma$ and $2n(1-\sqrt{1-h/n})\to h$, we may choose an integer $n>h$ with $2(1-\sqrt{1-h/n})<\gamma/n$.
By \Cref{lem small-gain}, there exists a positive rank-one operator $L$ such that
\begin{equation*}
	d(L)=1+h/n,\qquad\|L\|_1=\gamma/n.
\end{equation*}
With the notation of \Cref{thm main}, we have $\dim_{\mathcal P}(z_{\mathcal P}\mathcal H)=d(L)\geqslant 1$.
Thus, there exists an isometry $U\in\mathcal B(\mathcal H)$ with $U\mathcal H\subseteq z_{\mathcal P}\mathcal H$ such that
\begin{equation}\label{equ SU=Us}
	(s_j+\delta_{j1}L)U=Us_j
	\qquad(j\in\Lambda).
\end{equation}
Let $L_0=0$ and define $L_m=L+UL_{m-1}U^*$ for $1\leqslant m\leqslant n$.
Then each $L_m$ is positive and has rank at most $m$.
Since $U^*U=I$, we have
\begin{equation*}
	\|L_m\|_1=\operatorname{Tr}(L_m)
	=\operatorname{Tr}(L)+\operatorname{Tr}(L_{m-1})
	=\|L_{m-1}\|_1+\gamma/n.
\end{equation*}
For every $j\in\Lambda$, the intertwining relation \eqref{equ SU=Us} gives
\begin{equation*}
	(s_j+\delta_{j1}L_m)U
	=(s_j+\delta_{j1}L)U+\delta_{j1}UL_{m-1}
	=U(s_j+\delta_{j1}L_{m-1}).
\end{equation*}
Thus, $U\mathcal H$ reduces every $s_j+\delta_{j1}L_m$ and the unitary operator $U\colon\mathcal H\to U\mathcal H$ satisfies
\begin{equation*}
	(s_j+\delta_{j1}L_m)|_{U\mathcal H}U=U(s_j+\delta_{j1}L_{m-1})
	\qquad(j\in\Lambda).
\end{equation*}
Since $U^*$ vanishes on $(U\mathcal H)^\perp$, the recursion for $L_m$ gives
\begin{equation*}
	(s_j+\delta_{j1}L_m)|_{(U\mathcal H)^\perp}
	=(s_j+\delta_{j1}L)|_{(U\mathcal H)^\perp}
	\qquad(j\in\Lambda).
\end{equation*}
For $m=1$, \eqref{equ SU=Us} shows that $U\mathcal H$ contributes dimension $1$, and hence $(U\mathcal H)^\perp$ contributes $d(L)-1=h/n$.
Since the restricted family on $(U\mathcal H)^\perp$ is independent of $m$, additivity gives $d(L_m)=d(L_{m-1})+h/n$ for $1\leqslant m\leqslant n$.
It follows that
\begin{equation*}
	d(L_n)=1+h,\qquad\|L_n\|_1=\gamma.
\end{equation*}
We complete the proof by taking $K_1=L_n$.
\end{proof}

\begin{proof}[Proof of \Cref{thm multiplicity-range}]
If $\gamma=0$, then clearly $d=1$.
Suppose that $0<\gamma<\infty$.
By \Cref{lem multiplicity}, we have $1\leqslant d<1+\gamma$.
Conversely, every number in the interval $[1,1+\gamma)$ can be realized by \Cref{prop realization}.
\end{proof}

%\section*{Acknowledgments}

\end{document}